\documentclass[12pt,leqno]{article}
\usepackage{indentfirst}

\usepackage{amsmath}
\usepackage{amscd}
\usepackage{amssymb}
\usepackage{amsthm}

\def\zR{\ensuremath{\mathbb{R}}}

\newcommand{\Proof}{{\it Proof}}

\numberwithin{equation}{section}
\renewcommand{\theequation}{\thesection.\arabic{equation}}
\makeatletter
\newtheoremstyle{estiloobs}
  {6pt}
  {0pt}
  {\upshape}
  {}
  {}
  {}
  {.5em}
  {(\thmnumber{#2}) \thmname{\textbf{#1}}
   \thmnote{({#3})}}
\makeatother
\theoremstyle{estiloobs}
\newtheorem{remark}[equation]{Remark:}

\makeatletter
\newtheoremstyle{estilotheorem}
  {6pt}
  {0pt}
  {\slshape}
  {}
  {}
  {}
  {.5em}
  {(\thmnumber{#2}) \thmname{\textbf{#1}}
   \thmnote{({#3})}}
\makeatother
\theoremstyle{estilotheorem}

\newtheorem{proposition}[equation]{Proposition:}
\newtheorem{lemma}[equation]{Lemma:}

\makeatletter
\newtheoremstyle{estilolista}
  {0pt}
  {0pt}
  {\slshape}
  {}
  {}
  {}
  {.5em}
  {(\thmnumber{#2})\thmname{\textbf{#1}}
   \thmnote{({#3})}}
\makeatother
\theoremstyle{estilolista}

\newcommand{\pedro}{\theequation}

\begin{document}

\title{Boundedness of the fractional maximal operator on variable exponent Lebesgue spaces: a short proof}

\vskip 0.3 truecm

\author{Osvaldo Gorosito, Gladis Pradolini and Oscar Salinas\thanks{Research
supported by Consejo Nacional de Investigaciones Cient\'\i ficas y
T\'ecnicas de la Rep\'ublica Argentina and Universidad Nacional
del Litoral.\newline \indent Keywords and phrases: variable
spaces, maximal fractional operators.
}}

\date{\vspace{-1.5cm}}

\maketitle

\begin{abstract}
We give a simple proof of the boundedness of the fractional
maximal operator providing in this way an alternative approach to
the one given by C. Capone, D. Cruz Uribe and A. Fiorenza in
\cite{CCUF}.
\end{abstract}

\section{Introduction and main results} Given an open set $\Omega\subset\zR^n$
and a measurable function $p:\Omega \rightarrow [1,+\infty)$, the
variable exponent Lebesgue space $L^{p(.)}$ consists of all
measurable functions $f$ on $\Omega$ such that for some
$\lambda>0$, $
\int_{\Omega}\left(\frac{|f(x)|}{\lambda}\right)^{p(x)}\,
dx<\infty $ equipped with the norm
\begin{equation*}\|f\|_{p(.),\Omega}=\inf \{\lambda >0:
\int_{\Omega}\left(|f(x)|/\lambda\right)^{p(x)}\, dx\leq 1\}.
\end{equation*}

Given $\alpha$, $0<\alpha<n$, the fractional maximal operator
$M_{\alpha}$ is defined by
\begin{equation*}
M_{\alpha}f(x)=\sup_{B\ni x}\frac{1}{|B|^{1-\alpha/n}}\int_{B\cap
\Omega}|f(y)|\, dy,
\end{equation*}
where the supremum is taken over all balls $B$ which contain $x$.
When $\alpha=0$, $M_{0}=M$ is the Hardy-Littlewood maximal
operator.

In \cite{CCUF}, C. Capone, D. Cruz-Uribe and A. Fiorenza gave an
extension of the classical $L^p-L^q$ boundedness result for
$M_{\alpha}$ for $1<p<n/\alpha$ and $1/q=1/p-\alpha/n$ on the
variable Lebesgue context. The authors proved the following two
interesting pointwise inequalities relating both operators $M$ and
$M_{\alpha}$ which turn out to be crucial in proving the
boundedness result. Nevertheless, their proofs are not trivial at
all and they require additional lemmas. Moreover the hypotheses of
log-H\"{o}lder continuity on the exponent is needed.

\begin{proposition}{\rm ([CCUF])} Given an open set $\Omega\subset \zR^n$ and
$0<\alpha<n$, let $p:\Omega\rightarrow [1,\infty)$ be such that
$1<\inf_{x\in \Omega}p(x)\le p(x)\le \sup_{x\in
\Omega}p(x)<n/\alpha$ and such that
\begin{equation}\label{local}
|p(x)-p(y)|\leq \frac{C}{-\log |x-y|},\quad x,y\in \Omega, \quad
|x-y|<1/2.
\end{equation}
Let $q$ be such that $1/q(x)=1/p(x)-\alpha/n$. Then, for all $f\in
L^{p(.)}(\Omega)$ such that $\|f\|_{p(.),\Omega}\leq 1$ and such
that $f(x)\ge 1$ or $f(x)=0$, $x\in \Omega$,
\begin{equation*}
M_{\alpha}f(x)\leq CMf(x)^{p(x)/q(x)}.
\end{equation*}
\end{proposition}

\bigskip

\begin{proposition}{\rm ([CCUF])} Given an open set $\Omega\subset \zR^n$ and
$0<\alpha<n$, let $p:\Omega\rightarrow [1,\infty)$ be such that
$1<\inf_{x\in \Omega}p(x)\le p(x)\le \sup_{x\in
\Omega}p(x)<n/\alpha$ and such that
\begin{equation}\label{global}
|p(x)-p(y)|\leq \frac{C}{\log (e+ |x|)},\quad x,y\in \Omega, \quad
|y|\ge |x|.
\end{equation}
Let $q$ be such that $1/q(x)=1/p(x)-\alpha/n$. Then, for all $f\in
L^{p(.)}(\Omega)$ such that $\|f\|_{p(.),\Omega}\leq 1$ and such
that $0\leq f(x)< 1$, $x\in \Omega$,
\begin{equation*}
M_{\alpha}f(x)\leq CMf(x)^{p(x)/I_q(x)},
\end{equation*}
where $I_{q}(x)=\sup_{|y|\ge|x|}q(y)$
\end{proposition}

\medskip

Thus to obtain the boundedness result the authors have to split
the function $f$ into $f_1$ and $f_2$ properly and then apply the
propositions above to $M_{\alpha}(f_1)$ and $M_{\alpha}(f_2)$
respectively. Then the final result follows by applying the
continuity of $M$ proved in \cite{CUFN}.

\bigskip

The following elementary lemma is a successful substitute of the
propositions above. It should be also noticed that no conditions
of continuity on the exponent $p$ are required.

\begin{lemma}\label{nuestro}
Let $0<\alpha<n$ and $p$ be a function such that $1<\inf_{x\in
\Omega}p(x)\le p(x)\le \sup_{x\in \Omega}p(x)<n/\alpha$. Let $q$
be defined by $1/q(x)=1/p(x)-\alpha/n$. Then the following
inequality
\begin{equation*}
M_{\alpha}\left(f\right)(x)\leq
\left(M\left(|f|^{\frac{p(.)}{q(.)}\frac{n}{n-\alpha}}\right)(x)\right)^{1-
\alpha/n}\left(\int_{\Omega}|f(y)|^{p(y)}\, dy\right)^{\alpha/n}
\end{equation*}
holds for every function $f$.
\end{lemma}

\bigskip

\Proof\, : Let $Q\subset \zR^n$ be a cube containing $x$. Taking
into account that $p(y)/q(y)+\alpha p(y)/n=1$, by applying
H\"{o}lder's inequality we get
\begin{eqnarray*} \frac{1}{\mu(Q)^{1-\alpha
/n}}\int_{Q\cap\Omega}|f(y)|\, dy &=&\frac{1}{\mu(Q)^{1-\alpha
/n}}\int_{Q\cap\Omega}|f(y)|^{p(y)/q(y)}|f(y)|^{\alpha p(y)/n}\, dy\\
&\leq&
\left(\frac{1}{|Q|}\int_{Q\cap\Omega}|f(y)|^{\frac{p(y)}{q(y)}\frac{n}{n-\alpha}}\,
dy\right)^{1-
\alpha/n}\left(\int_{\Omega}|f(y)|^{p(y)}\, dy\right)^{\alpha/n}\\
&\leq&
\left(M\left(|f|^{\frac{p(.)}{q(.)}\frac{n}{n-\alpha}}\right)(x)\right)^{1-
\alpha/n}\left(\int_{\Omega}|f(y)|^{p(y)}\, dy\right)^{\alpha/n}.\\
\end{eqnarray*}
Thus the desired inequality follows inmediately.$\square$

\bigskip

A straightforward application of lemma $\ref{nuestro}$ allows us
to obtain the boundedness of $M_{\alpha}$ in the variable context.
In fact, since $f\in L^{p(x)}$ implies that
$|f|^{\frac{p(x)}{q(x)}\frac{n}{n-\alpha}}\in L^{q(x)(1-
\alpha/n)}$ the result follows from the boundedness of $M$.
Particularly if $p$ satisfies ($\ref{local}$) and ($\ref{global}$)
this was proved in \cite{CUFN}.

\medskip

\begin{remark}
A weighted pointwise inequality in a more general context was
proved in \cite{GPS}. In that paper the authors took advantage of
this result to obtain weighted results for the boundedness of the
fractional maximal operator in the variable context with a
non-necessary doubling measure.
\end{remark}

\end{document}